\documentclass[11pt]{article}

\usepackage{amsfonts}
\usepackage{amsmath}
\usepackage{amsthm}
\usepackage{amssymb}

\usepackage{cancel}
\usepackage{comment}

\theoremstyle{plain} 
\newtheorem*{theorem}{Theorem}

\newtheorem*{lemma*}{Lemma}
\newtheorem*{claim*}{Claim}
\newtheorem*{theorem*}{Theorem}
\newtheorem*{proposition*}{Proposition}
\newtheorem*{problemstatement*}{Problem}
\newtheorem*{corollary*}{Corollary}

\theoremstyle{definition} 
\newtheorem{definition}{Definition}

\newtheorem{example}[definition]{Example}

\newcommand{\ignorethis}[1]{}

\usepackage{graphicx}
\usepackage{ifpdf}
\ifpdf
\fi
\usepackage{cancel}
\usepackage{hyperref}
\newcommand{\inner}[2]{\left\langle #1 , #2 \right\rangle}
\usepackage{url}
\usepackage[numbers,sort]{natbib}
\usepackage{paralist}
\usepackage{enumitem}
\usepackage{verbatim} 

\begin{document}

\title{On complex H-type Lie algebras}
\author{Nathaniel Eldredge
  \footnote{
    School of Mathematical Sciences,
    University of Northern Colorado,
    501 20th St., Campus Box 122,
    Greeley, \textsc{CO}, 80639,
    \textsc{USA}.
    Email: \texttt{neldredge@unco.edu}.
    This work was supported by a grant from the Simons Foundation (\#355659, Nathaniel Eldredge).
 }
}
\maketitle

\begin{abstract}

    H-type Lie algebras were introduced by Kaplan as a class of real Lie
algebras generalizing the familiar Heisenberg Lie algebra
$\mathfrak{h}^3$.  The H-type property depends on a choice of inner
product on the Lie algebra $\mathfrak{g}$.  Among the H-type Lie
algebras are the complex Heisenberg Lie algebras
$\mathfrak{h}^{2n+1}_{\mathbb{C}}$, for which the standard Euclidean
inner product not only satisfies the H-type condition, but is also
compatible with the complex structure, in that it is Hermitian. We
show that, up to isometric isomorphism, these are the only complex Lie
algebras with an inner product satisfying both conditions.  In other
words, the family $\mathfrak{h}^{2n+1}_{\mathbb{C}}$ comprises all of
the complex H-type Lie algebras.

\end{abstract}

\section{Introduction}

Since their introduction by Kaplan \cite{kaplan80}, H-type Lie
algebras, and their corresponding nilpotent Lie groups, have attracted
interest as a natural generalization of the classical real Heisenberg
Lie algebra $\mathfrak{h}^3$ of dimension 3 and the corresponding real
Heisenberg group $\mathbb{H}^3$.  The Heisenberg group is a motivating
example in many areas of mathematics, and in many cases, facts about
the Heisenberg group carry over into the H-type setting.  For
instance, H-type groups carry a natural structure as sub-Riemannian
manifolds, and the analysis of their sub-Laplacians has attracted
considerable interest.  As a sampling, we mention
\cite{inglis-kontis-zegarlinski,eldredge-gradient,reimann-rigidity,li-gradient-h-type,balogh-tyson,he-yin-hardy}.

The H-type condition for a (real) Lie algebra $\mathfrak{g}$ is
dependent on a choice of inner product $\langle \cdot, \cdot \rangle$
(i.e. a positive definite, symmetric, bilinear form) on
$\mathfrak{g}$, so it is really a property of the pair $(\mathfrak{g},
\langle \cdot, \cdot \rangle)$.  For example, in $\mathfrak{h}^3$, the natural
Euclidean inner product will do.

In this note, we focus on the complex Heisenberg (or Heisenberg--Weyl)
Lie algebras $\mathfrak{h}^{2n+1}_{\mathbb{C}}$, which, when
considered as real Lie algebras and equipped with their natural
Euclidean inner products, are likewise H-type.  But these Lie
algebras also carry a complex structure, and the Euclidean inner
product is Hermitian with respect to this structure, which is a
natural compatibility condition.  As such, analysis on the complex
Heisenberg groups $\mathbb{H}^{2n+1}_{\mathbb{C}}$ can take advantage
of all the tools of complex geometry, together with the many results
for H-type groups mentioned above.  However, the purpose of this note
is to show that this harmonious relationship between these structures is
essentially unique to these specific Lie algebras (and their
respective Lie groups).

As an application, we refer to
\cite{eldredge-gross-saloff-coste-dila}, in which we studied a
property known as strong hypercontractivity for the hypoelliptic heat
kernel on a stratified complex Lie group.  An essential hypothesis for
this result was that the heat kernel should satisfy a logarithmic
Sobolev inequality.  For most Lie groups, it remains an open problem
to determine whether this inequality holds, but it follows from the
results of \cite{eldredge-gradient,li-gradient-h-type} that the
inequality holds in every H-type Lie group.  Thus, the strong
hypercontractivity theorem proved in
\cite{eldredge-gross-saloff-coste-dila} holds in particular for any
complex Lie group which, when considered as a real Lie group, is also
H-type.  The result of the present note implies that these Lie groups
are precisely the family $\mathbb{H}^{2n+1}_{\mathbb{C}}$.  As this is
a relatively limited class of examples, we see this as further
motivation to try to extend the logarithmic Sobolev inequality beyond
the H-type case.

\section{Definitions and examples}

We begin by recalling the definition of an H-type Lie algebra, as
formulated in \cite[Definition 18.1.1]{blu-book}. (Kaplan's original
definition \cite{kaplan80} is equivalent, but slightly less
convenient for our purposes.)  Let $\mathfrak{g}$ be a real
finite-dimensional Lie algebra equipped with an inner product
$\inner{\cdot}{\cdot} : \mathfrak{g} \times \mathfrak{g} \to
\mathbb{R}$.  Let $\mathfrak{z}$ be the center of $\mathfrak{g}$, and
let $\mathfrak{v} = \mathfrak{z}^\perp$.  For $z \in \mathfrak{z}$ and
$u \in \mathfrak{v}$, define $J_z u$ as the unique element of
$\mathfrak{v}$ satisfying
\begin{equation}\label{Jdef}
  \inner{J_z u}{v} = \inner{z}{[u,v]} \quad \text{for all $v \in \mathfrak{v}$}.
\end{equation}
It is clear that each $J_z : \mathfrak{v} \to \mathfrak{v}$ is a linear
map, and moreover $z \mapsto J_z$ is linear in $z$.

\begin{definition}\label{Htype-def}
We say that $(\mathfrak{g}, \inner{\cdot}{\cdot})$ is \textbf{H-type}
if the following two conditions hold:
\begin{enumerate}
\item $[\mathfrak{v}, \mathfrak{v}] \subset \mathfrak{z}$
\item For each $z \in \mathfrak{z}$ with $\|z\|=1$, $J_z :
  \mathfrak{v} \to \mathfrak{v}$ is an isometry with respect to
  $\inner{\cdot}{\cdot}$.
\end{enumerate}
\end{definition}
We observe that an H-type Lie algebra is necessarily nilpotent of step
2.  A simply-connected Lie group is said to be H-type if its Lie
algebra is H-type in the above sense.

Now suppose that $\mathfrak{g}$ is a complex Lie algebra, whose
complex structure we denote by $i$.  If we wish to equip
$\mathfrak{g}$ with a real inner product, it is natural to demand some
compatibility with the complex structure.  Specifically, we would like
the inner product to be \textbf{Hermitian}, i.e., for $x,y \in
\mathfrak{g}$ we have $\inner{ix}{iy} = \inner{x}{y}$.  We may then
define $J$ in terms of this inner product by (\ref{Jdef}).  We
observe for later use that, as a consequence of the Hermitian property
of the inner product, we have for $\alpha, \beta \in \mathbb{C}$ and
$u,z \in \mathfrak{g}$,
\begin{equation}\label{complex-J}
  J_{\alpha z}(\beta u) = \alpha \bar{\beta} J_z u.
\end{equation}
That is, $J_z u$ is complex linear in $z$ and conjugate linear in
$u$.

The question of interest in this note is when both of the above
properties hold, motivating the following definition.

\begin{definition}
  A \textbf{complex H-type Lie algebra} is a pair $(\mathfrak{g},
  \langle \cdot, \cdot \rangle)$, where $\mathfrak{g}$ is a complex Lie
  algebra and $\langle \cdot, \cdot \rangle$ is an inner product on
  $\mathfrak{g}$, such that the following two conditions hold:
  \begin{itemize}
  \item The inner product $\langle \cdot, \cdot \rangle$ is Hermitian
    with respect to the complex structure of $\mathfrak{g}$.
  \item Forgetting the complex structure on $\mathfrak{g}$, the pair
    $(\mathfrak{g}, \langle \cdot, \cdot \rangle)$ is H-type in the
    sense of  Definition \ref{Htype-def}.
  \end{itemize}
\end{definition}

We can likewise define a \textbf{complex H-type Lie group} as a
connected and simply connected complex Lie group $G$ equipped with a
Hermitian left-invariant Riemannian metric $g$ which, when viewed as
an inner product on the Lie algebra of $G$, satisfies the above
conditions.

\begin{example}\label{complex-heis}
  The \textbf{complex Heisenberg Lie algebra} of complex
  dimension $2n+1$ is the complex Lie algebra
  $\mathfrak{h}^{2n+1}_{\mathbb{C}}$ generated (over $\mathbb{C}$) by
  the basis of the $2n+1$ vectors $\{x_1, y_1, \dots, x_n, y_n, z\}$, with the bracket defined
  by $[x_k, y_k]=z$, and for $j \ne k$, $[x_j, y_k] = [x_j, z] = [y_j,
    z] = 0$.  We may equip $\mathfrak{h}^{2n+1}_{\mathbb{C}}$ with the
  real inner product $\inner{\cdot}{\cdot}$ that makes all of $x_k,
  ix_k, y_k, iy_k, z, iz$ orthonormal; it is clear that this inner
  product is Hermitian.  The center $\mathfrak{z}$ of
  $\mathfrak{h}^{2n+1}_{\mathbb{C}}$ is spanned (over $\mathbb{C}$) by
  $z$, so we clearly have $[\mathfrak{v}, \mathfrak{v}] =
  \mathfrak{z}$.  Defining $J$ as above, it is easy to compute
  \begin{equation*}
    J_z x_k = y_k \quad J_z y_k = -x_k \quad J_z ix_k = -iy_k
    \quad J_z iy_k = ix_k
  \end{equation*}
  so that $J_z$ is an isometry.  Moreover, every element $w \in
  \mathfrak{z}$ is of the form $w = \alpha z$ for some $\alpha \in \mathbb{C}$, and $\|w\| = |\alpha|$,
  so using (\ref{complex-J}) we see that $J_w$ is an isometry whenever
  $\|w\| = 1$.  Thus $(\mathfrak{h}^{2n+1}_{\mathbb{C}},
  \inner{\cdot}{\cdot})$ is a  complex H-type Lie algebra.
\end{example}

Of course, the complex Heisenberg Lie algebras are a very special
family within the far larger class of all complex Lie algebras.
Likewise, the class of H-type Lie algebras, although fairly
restrictive, is still much broader than this specific family.  For
instance, there exist H-type Lie algebras having centers of any given
real dimension \cite{kaplan80}, while the complex Heisenberg Lie
algebras all have centers of real dimension $2$.

Nevertheless, we shall now prove that the complex Heisenberg Lie
algebras are, up to isometric
isomorphism, the only complex H-type Lie algebras.

\section{Main result}

\begin{theorem}
  Let $(\mathfrak{g}, \inner{\cdot}{\cdot})$ be a complex H-type Lie
  algebra as defined above.  Then for some $n$, $(\mathfrak{g},
  \inner{\cdot}{\cdot})$ is isometrically isomorphic to
  $\mathfrak{h}^{2n+1}_{\mathbb{C}}$ with its standard Hermitian inner
  product.
\end{theorem}

In particular, complex H-type Lie algebras are completely classified
by their dimension.  We also immediately obtain the analogous
classification of complex H-type Lie groups.

\begin{proof}
  Suppose $(\mathfrak{g}, \inner{\cdot}{\cdot})$ is complex H-type,
  and let $\mathfrak{v}$, $\mathfrak{z}$ and $J$ be defined as above.

  We recall the well-known Clifford algebra identity for H-type Lie
  algebras:
  \begin{equation}\label{clifford}
    J_z J_w + J_w J_z = -2 \inner{z}{w} I, \quad z,w \in \mathfrak{z}.
  \end{equation}
  To prove this, first consider the case when $w=z$ and $\|z\|=1$.  Then
  for any $u,v \in \mathfrak{v}$, we have
  \begin{align*}
    \inner{J_z^2 u}{v} = \inner{z}{[J_z u, v]} = -\inner{z}{[v, J_z u]}
  = -\inner{J_z v}{J_z u} = -\inner{v}{u}
  \end{align*}
  since $J_z$ is an isometry.  So $J_z^2 = -I$.  The general case
  follows by scaling and polarization.
  
  We begin by showing that $\mathfrak{z}$ must have complex dimension
  1.  If not, then we can find $z,w \in \mathfrak{z}$ with $\|z\| =
  \|w\| = 1$ and $\inner{z}{w} = \inner{iz}{w} = 0$.  Then by
  (\ref{clifford}) and (\ref{complex-J}) we have
  \begin{align*}
    0 &= -2 \inner{z}{w} I = J_z J_w + J_w J_z \\
    0 &= -2 \inner{iz}{w} I = J_{iz} J_w + J_w J_{iz} = i J_z J_w +
    J_w i J_z = i(J_z J_w - J_w J_z).
  \end{align*}
  Thus $J_w J_z = J_z J_w = 0$, contradicting the requirement that
  $J_z, J_w$ be isometries.

  Therefore, $\mathfrak{z}$ is the complex span of a single unit
  vector $z$.  We recursively construct an orthonormal basis for
  $\mathfrak{v}$ over $\mathbb{R}$, of the form $\{x_k, ix_k,
  y_k, iy_k : k=1,\dots,n\}$.  Suppose $\{x_k, ix_k, y_k, iy_k : k =
  1, \dots, m-1\}$ have been constructed and do not span
  $\mathfrak{v}$.  Let $x_m$ be any unit vector orthogonal to all of
  $x_k, ix_k, y_k, iy_k$ for $k=1,\dots,m$.  Then set $y_m = J_z x_m$.
  We have $\|y_m\| = 1$, and a few straightforward computations verify
  that $\{x_k, ix_k, y_k, iy_k : k = 1,\dots, m\}$ are now orthogonal.
  When the process terminates, we have the desired orthonormal basis.

  To compute brackets, for $j \ne k$ we have
  \begin{align*}
    \inner{z}{[x_k, y_k]} &= \inner{J_z x_k}{y_k} = \inner{y_k}{y_k}
    = 1 \\
    \inner{z}{[x_k, x_j]} &= \inner{J_z x_k}{x_j} = \inner{y_k}{x_j}
    = 0 \\
    \inner{z}{[y_k,y_j]} &= \inner{J_z y_k}{y_j} = \inner{J_z y_k}{J_z
      x_j} = \inner{y_k}{x_j} = 0 \\
    \inner{z}{[x_k, y_j]} &= \inner{J_z x_k}{y_j} =
    \inner{y_k}{y_j} = 0.
  \end{align*}
  Similar computations show that if $z$ is replaced by $iz$, all of
  the above expressions vanish.  Each bracket is in $\mathfrak{z}$ and
  hence a complex scalar multiple of $z$, so we have
  \begin{equation*}
    [x_k, y_k] = z, \quad [x_k, x_j] = [y_k, y_j] = [x_k, y_j] = 0.
  \end{equation*}
  The corresponding brackets for $ix_k, iy_k$, etc, follow from the
  complex bilinearity of the bracket.  These are precisely the
  same relations as for the complex Heisenberg Lie algebra
  $\mathfrak{h}^{2n+1}_{\mathbb{C}}$, and the basis is orthonormal,
  just as for the standard inner product on
  $\mathfrak{h}^{2n+1}_{\mathbb{C}}$.
  Therefore, the unique complex
  linear map $\mathfrak{g} \to \mathfrak{h}^{2n+1}_{\mathbb{C}}$
  sending $x_1, y_1, \dots, x_n, y_n, z \in \mathfrak{g}$ to the standard basis for
  $\mathfrak{h}^{2n+1}_{\mathbb{C}}$ (described in Example
  \ref{complex-heis}) is an isometric isomorphism of complex Lie
  algebras.
\end{proof}


\emph{Acknowledgments.}
The author would like to thank Maria
Gordina, Leonard Gross, Martin Moskowitz, and Laurent Saloff-Coste for
helpful suggestions.

\bibliographystyle{plainnat}
\bibliography{allpapers}

\def\cprime{$'$} \def\cprime{$'$} \def\cprime{$'$} \def\cprime{$'$}
\begin{thebibliography}{9}
\providecommand{\natexlab}[1]{#1}
\providecommand{\url}[1]{\texttt{#1}}
\expandafter\ifx\csname urlstyle\endcsname\relax
  \providecommand{\doi}[1]{doi: #1}\else
  \providecommand{\doi}{doi: \begingroup \urlstyle{rm}\Url}\fi

\bibitem[Balogh and Tyson(2002)]{balogh-tyson}
Zolt{\'a}n~M. Balogh and Jeremy~T. Tyson.
\newblock Polar coordinates in {C}arnot groups.
\newblock \emph{Math. Z.}, 241\penalty0 (4):\penalty0 697--730, 2002.
\newblock ISSN 0025-5874.
\newblock \doi{10.1007/s00209-002-0441-7}.
\newblock URL \url{http://dx.doi.org/10.1007/s00209-002-0441-7}.

\bibitem[Bonfiglioli et~al.(2007)Bonfiglioli, Lanconelli, and
  Uguzzoni]{blu-book}
A.~Bonfiglioli, E.~Lanconelli, and F.~Uguzzoni.
\newblock \emph{Stratified {L}ie groups and potential theory for their
  sub-{L}aplacians}.
\newblock Springer Monographs in Mathematics. Springer, Berlin, 2007.
\newblock ISBN 978-3-540-71896-3; 3-540-71896-6.

\bibitem[Eldredge(2010)]{eldredge-gradient}
Nathaniel Eldredge.
\newblock Gradient estimates for the subelliptic heat kernel on {$H$}-type
  groups.
\newblock \emph{J. Funct. Anal.}, 258\penalty0 (2):\penalty0 504--533, 2010.
\newblock ISSN 0022-1236.
\newblock \doi{10.1016/j.jfa.2009.08.012}.
\newblock URL \url{http://dx.doi.org/10.1016/j.jfa.2009.08.012}.
\newblock arXiv:0904.1781.

\bibitem[Eldredge et~al.(2015)Eldredge, Gross, and
  Saloff-Coste]{eldredge-gross-saloff-coste-dila}
Nathaniel Eldredge, Leonard Gross, and Laurent Saloff-Coste.
\newblock Strong hypercontractivity and logarithmic {S}obolev inequalities on
  stratified complex {L}ie groups.
\newblock Preprint. arXiv:1510.05151, 2015.
\newblock URL \url{http://arxiv.org/abs/1510.05151}.

\bibitem[He and Yin(2016)]{he-yin-hardy}
Jianxun He and Mingkai Yin.
\newblock {$L^p$} {H}ardy type inequality in the half-space on the {H}-type
  group.
\newblock \emph{J. Inequal. Appl.}, pages 2016:129, 10, 2016.
\newblock ISSN 1029-242X.
\newblock \doi{10.1186/s13660-016-1070-8}.
\newblock URL \url{http://dx.doi.org/10.1186/s13660-016-1070-8}.

\bibitem[Hu and Li(2010)]{li-gradient-h-type}
Jun-Qi Hu and Hong-Quan Li.
\newblock Gradient estimates for the heat semigroup on {H}-type groups.
\newblock \emph{Potential Anal.}, 33\penalty0 (4):\penalty0 355--386, 2010.
\newblock ISSN 0926-2601.
\newblock \doi{10.1007/s11118-010-9173-1}.
\newblock URL \url{http://dx.doi.org/10.1007/s11118-010-9173-1}.

\bibitem[Inglis et~al.(2011)Inglis, Kontis, and
  Zegarli{\'n}ski]{inglis-kontis-zegarlinski}
J.~Inglis, V.~Kontis, and B.~Zegarli{\'n}ski.
\newblock From {$U$}-bounds to isoperimetry with applications to {H}-type
  groups.
\newblock \emph{J. Funct. Anal.}, 260\penalty0 (1):\penalty0 76--116, 2011.
\newblock ISSN 0022-1236.
\newblock \doi{10.1016/j.jfa.2010.08.003}.
\newblock URL \url{http://dx.doi.org/10.1016/j.jfa.2010.08.003}.

\bibitem[Kaplan(1980)]{kaplan80}
Aroldo Kaplan.
\newblock Fundamental solutions for a class of hypoelliptic {PDE} generated by
  composition of quadratic forms.
\newblock \emph{Trans. Amer. Math. Soc.}, 258\penalty0 (1):\penalty0 147--153,
  1980.
\newblock ISSN 0002-9947.

\bibitem[Reimann(2001)]{reimann-rigidity}
Hans~Martin Reimann.
\newblock Rigidity of {$H$}-type groups.
\newblock \emph{Math. Z.}, 237\penalty0 (4):\penalty0 697--725, 2001.
\newblock ISSN 0025-5874.
\newblock \doi{10.1007/PL00004887}.
\newblock URL \url{http://dx.doi.org/10.1007/PL00004887}.

\end{thebibliography}

\end{document}